\documentstyle[12pt, fullpage, amsfonts]{amsart}
\begin{document}
\thispagestyle{empty}
\vfuzz=2pt

\newtheorem{theorem}{Theorem}
\newtheorem{lemma}{Lemma}
\newtheorem{corollary}{Corollary}
\newtheorem{proposition}[lemma]{Proposition}
\newtheorem{fact}{Fact}

\newenvironment{pflike}[1]{\medskip\noindent{\bf #1}}{\medskip}
\newenvironment{proof}{\noindent{\bf Proof:}}{\qed\medskip}

\newcommand{\R}{{\mathbb R}}

\newcommand{\A}{{\cal A}}
\newcommand{\C}{{\cal C}}
\newcommand{\M}{{\cal M}}
\renewcommand{\P}{{\cal P}}
\renewcommand{\S}{{\cal S}}

\newcommand{\Int}{\mathop{\rm Int}}
\newcommand{\vol}{\mathop{\rm vol}}
\newcommand{\diam}{\mathop{\rm diam}}

\newcommand{\ep}{\varepsilon}
\newcommand{\onrn}{O(n)\times\R^n}
\renewcommand{\j}{_j} % {^{(j)}}

\title{Compactness Theorems for Geometric Packings}
\author{Greg Martin}
\address{Department of Mathematics\\University of Toronto\\Canada M5S 3G3}
\email{gerg@@math.toronto.edu}
\begin{abstract}
Moser asked whether the collection of rectangles of
dimensions $1\times\frac12$, $\frac12\times\frac13$,
$\frac13\times\frac14$, \dots, whose total area equals 1, can be
packed into the unit square without overlap, and whether the
collection of squares of side lengths $\frac12$, $\frac13$, $\frac14$,
\dots can be packed without overlap into a rectangle of area
$\frac{\pi^2}6-1$. Computational investigations have been made into
packing these collections into squares of side length $1+\ep$ and
rectangles of area $\frac{\pi^2}6-1+\ep$, respectively, and one can
consider the apparently weaker question of whether such packings are
possible for every positive number~$\ep$. In this paper we establish a
general theorem on sequences of geometrical packings that implies, in
particular, that the ``for every~$\ep$'' versions of these two problems
are actually equivalent to the original tiling problems.
\end{abstract}
\maketitle

\section{Introduction}

\noindent Given a collection $\A=\{A_1,A_2,\dots\}$ of subsets of $\R^n$, a {\it
packing\/} of $\A$ into another set $C\subset\R^n$ is a way of fitting
each of the sets $A_i$ inside $C$ without overlap. By a positioning of a set
$A_i$ we mean the image of $A_i$ under a rigid motion of
$\R^n$, i.e., some combination of translations, rotations, and reflections. To
avoid ambiguity about points on the boundaries of the $A_i$, we say more
precisely that these positionings of the $A_i$ must be contained inside $C$ and
that their interiors must be pairwise disjoint. One can also speak of {\it
oriented packings\/}, where the sets $A_i$ may be translated and rotated but not
reflected, and also {\it translated packings\/}, where the $A_i$ may be
translated but neither rotated nor reflected. We also refer to a translated
packing as a {\it parallel packing\/}, particularly when each set $A_i$ is a
brick (a product $[x_1,y_1]\times\dots\times[x_n,y_n]$ of closed
intervals). If the union of the repositioned sets $A_i$ is all of $C$, we call
the packing a {\it tiling\/} of $C$.

It is often difficult to determine whether a particular collection
$\A$ can be packed into some target set $C$. One representative
example is the collection $\A=\{A_1,A_2,\dots\}$ where each $A_i$ is a
rectangle of dimensions $\frac1i\times\frac1{i+1}$. Since the total
area of these rectangles is~1, it is conceivable that $\A$ can tile a unit square
(generally or even with a parallel tiling); but this
problem, first posed by Moser (see~\cite{meirmoser} and~\cite[Section
D5]{UPIG}), is unsolved. One can instead ask the apparently weaker question
of whether for every positive number $\ep$, the collection $\A$ can be packed
inside a square of side length $1+\ep$ (see for
example~\cite{ball}). A similar situation holds with the collection
$\A=\{S_2,S_3,\dots\}$ where each $S_i$ is a square of side length
$\frac1i$. Conceivably this collection will tile a rectangle of area
${\pi^2\over6}-1$ (and perhaps even one with dimensions
$({\pi^2\over6}-1)\times1$), but it is even unknown whether for every
positive number $\ep$ the collection $\A$ can be packed into
rectangles with area ${\pi^2\over6}-1+\ep$. For both these problems,
results of Paulhus~\cite{paulhus} shows that $\ep$ can at least be
taken smaller than $10^{-9}$.

The purpose of this paper is to show that the weaker ``for every~$\ep$'' versions
of these two packing problems are actually equivalent to the stronger tiling
versions. Our methods apply in a somewhat more general setting,
and we state the following two theorems as representative of what can be
deduced. %\vadjust{\eject}
For the first theorem, we use the notation $\lambda C = \{\lambda
y\colon y\in C\}$ for the homothetic expansion/dilation (or simply {\it
homothet\/}) of $C$ by the constant factor $\lambda>0$.

\begin{theorem} Let $\A$ be a collection of subsets of $\R^n$, and let $C$ be a
compact subset of $\R^n$. If for every $\ep>0$ there exists a packing of $\A$
into the homothet $(1+\ep)C$, then there exists a packing of $\A$
into $C$ itself. In particular, if there exist packings of $\A$ into closed
balls of radius $R+\ep$ for every $\ep>0$, then there exists a packing of $\A$
into a closed ball of radius $R$. These statements remain true if ``packing'' is
replaced by ``oriented packing'' or ``translated packing''.
\label{homothetthm}
\end{theorem}

\noindent We remark that the collection $\A$ may have any cardinality. Of course,
the hypothesis that the target set $C$ be compact is equivalent to $C$ being
both closed and bounded; both of these conditions on $C$ are necessary. There
are obvious counterexamples if $C$ is not required to be closed---for example, we
can take $C$ to be the open unit disk in $\R^2$ and $\A$ to be the collection
consisting solely of $\bar C$, the closure of $C$. The theorem also fails if $C$
is closed but not bounded: for example, we can again take $\A$ to
consist solely of the closed unit disk in $\R^2$, and $C$ to be the
the closed region $\{(x,y)\colon 1\le x,\, |y|\le1-1/x\}$.

\begin{theorem} Let $\A$ be a collection of subsets of $\R^n$. If there exist
packings of $\A$ into bricks of volume $V+\ep$ for every $\ep>0$, then there
exists a packing of $\A$ into a brick of volume~$V$. In fact a stronger
statement is true: let $\{B_1,B_2,\dots\}$ be a sequence of bricks in $\R^n$,
with the dimensions of the $j$th brick $B_j$ being $b_{j1}\times\dots\times
b_{jn}$. Set $V = \inf_j\{\vol B_j\}$, and assume that $\vol B_j > V$ for every
$j$. Suppose that there exists a packing of $\A$ into each brick $B_j$. Then
there exists a packing of $\A$ into some brick $B$ with dimensions
$b_1\times\dots\times b_n$, satisfying $\vol B = V$ and $b_m \le \limsup_j
\{b_{jm}\}$ for each $1\le m\le n$. These statements remain true if ``packing''
is replaced by ``oriented packing'' or ``translated packing''.
\label{brickthm}
\end{theorem}

The equivalence of the weak and strong versions of the two packing problems
mentioned in the introductory remarks follow as immediate corollaries of
Theorem~\ref{brickthm}:

\begin{corollary} Let $\A$ be the collection of rectangles of dimensions
$1\times{1\over2}$,
${1\over2}\times{1\over3}$, ${1\over3}\times{1\over4}$,
${1\over4}\times{1\over5}$, \dots. Suppose that for every $\ep>0$, the
collection $\A$ can be packed into a square of area $1+\ep$. Then $\A$ tiles a
square of area 1. If the given packings are parallel packings, then
$\A$ parallel-tiles a square of area 1.
\end{corollary}

\begin{corollary} Let $\A$ be the collection of squares of side lengths
${1\over2}$,
${1\over3}$, ${1\over4}$, \dots. Suppose that for every $\ep>0$, the collection
$\A$ can be packed into a rectangle of area
${\pi^2\over6}-1+\ep$. Then $\A$ tiles a rectangle of area ${\pi^2\over6}-1$. If
the given packings are into rectangles of height 1, then $\A$ tiles a rectangle
of dimensions $1\times({\pi^2\over6}-1)$. In either case, if the given packings
are parallel packings, then $\A$ parallel-tiles the resulting rectangle of area
${\pi^2\over6}-1$.
\end{corollary}

The aforementioned work of Paulhus~\cite{paulhus} makes a convincing argument
that the ``for every~$\ep$'' versions of these two packing questions have
affirmative answers (since obstacles to finding rectangle tilings
generally arise from the largest rectangles). In light of Corollaries~1 and ~2,
it therefore seems likely that tilings (indeed, parallel tilings) do exist in
both cases.
 
As can be inferred from the title of this paper, the methods used to establish
Theorems~1 and ~2 are topological in nature. The intuitive idea is to convert a
sequence of packings of the collection $\A$ in the hypothesized sets into a
``limiting packing'' of $\A$ into the desired target set. To this end, we will
show how the set of packings of $\A$ can be naturally regarded as a
topological space, and then use a compactness argument to show the existence of
a ``limiting packing'' of some sort; it then remains to show that this packing
is a valid packing into the type of set required by Theorem~1 or~2.

In Section~2 we set the notation to be used throughout this paper and exhibit
simple properties of the defined objects that follow easily from elementary
point-set topology. Section~3 contains the proofs of Theorems~1 and ~2, modulo
an important proposition whose proof will be deferred until Section~4 in order
to clarify the issues involved in the proofs of the theorems themselves. In
Section~5 we remark on some modified versions of Theorems~1 and~2 that
can be proved using these methods, without going into the details of the proofs.

\section{Notation and Basic Topological Facts}

\noindent The methods that we use are valid for collections $\A$ of
subsets of $\R^n$ of any cardinality, but for the sake of
notational simplicity we work under the assumption that our collection
$\A = \{A_1,A_2,\dots\}$ is countably infinite. In addition, we argue
throughout with the understanding that we are allowing translations,
rotations, and reflections and thus permitting the most general kinds
of packings; at the beginning of Section~5 we will explain how our
arguments extend to the more restrictive classes of oriented packings
and parallel packings.

For any subset $C$ of $\R^n$, we denote by $\P(\A,C)$ the set (possibly empty
{\it a priori\/}) of all packings of $\A$ into $C$. We mention at the outset
that translated copies of the target space $C$ are equivalent to each other for
the purposes of deciding whether there exists a packing of $\A$ into
$C$---indeed, there is a natural bijection between the set of packings of
$\A$ into $C$ and the set of packings of $\A$ into some translated copy of $C$.
Similarly, we may modify the collection $\A$ by replacing each set $A_i$ by
any translated copy of $A_i$, and still retain in essence the same set
$\P(\A,C)$. For instance, it will often be convenient for us to assume that each
set $A_i$ contains the origin in $\R^n$. We also note that if $C$ is a subset of
$D$ then certainly $\P(\A,C)
\subset
\P(\A,D)$.

Let $O(n)$ denote the $n$-dimensional orthogonal group, i.e., the set of all
$n\times n$ matrices $\theta$ with real entries such that $\theta^{-1} =
\theta^T$. Every rigid motion of $\R^n$ can be identified with an element of the
product space $\onrn$ as follows: if $\sigma=(\theta,\xi)$ is an element of
$\onrn$, then $\sigma$ acts on a point $x$ of $\R^n$ by the rule $\sigma(x) =
\xi + \theta x$. (Throughout this paper we will maintain the notational
conventions that elements of $\onrn$ will be denoted by $\sigma$ or $\tau$, and
that $\theta$ and $\xi$ will denote the $O(n)$- and $\R^n$-components,
respectively, when it is necessary to refer to these components separately.)
Certainly these rigid motions $\sigma$ act on subsets $A$ of $\R^n$ as well, and
we will write $\sigma(A) = \{\xi+\theta x\colon x\in A\}$ for the image. Any
positioning of the set $A$ in $\R^n$, using translations, rotations, and/or
reflections, can be realized as $\sigma(A)$ for some element $\sigma$ of $\onrn$.

Define the topological space $\M(\R^n)$ to be the product
space $(\onrn)^\infty$, and for any subset $D$ of $\R^n$ define the subspace
$\M(D) = (O(n)\times D)^\infty$ of $\M(\R^n)$. Since every positioning of a set
$A$ in $\R^n$ corresponds uniquely to an element $\sigma$ of $\onrn$, the space
$\M(\R^n)$ parametrizes all possible positionings of the collection $\A$ in
$\R^n$, and certain positionings among these will correspond to packings of $\A$
into a target set $C$. More precisely, if $\Int A$ denotes the interior of $A$,
we can write
\begin{equation}
\begin{split}
\P(\A,C) = \big\{ S = \{ \sigma_i \} \in \M(\R^n) \colon &\forall i,\,
\sigma_i(A_i) \subset C; \\
&\forall i\ne j,\, \Int(\sigma_i(A_i)) \cap \Int(\sigma_j(A_j)) = \emptyset
\big\}.
\end{split}
\label{pacspace}
\end{equation} \goodbreak\noindent
(In general we will let $S$ and $T$ denote elements of $\M(\R^n)$ or of its
subsets.) As a result, the set $\P(\A,C)$ can be given the subspace topology
induced by the product topology on $\M(\R^n)$. The key to the proof of
Theorem~\ref{homothetthm} is to exploit this topological structure on $\M(\R^n)$
to show that $\P(\A,C)$ is a nonempty subspace under the stated hypotheses, and
the proof of Theorem~\ref{brickthm} proceeds similarly after a suitable brick
$B$ is chosen as the ultimate target set.

We now exhibit several facts, which follow from the definitions of the above
notation together with elementary point-set topology, that will be useful to us
later. As a final piece of notation, let
\begin{equation*}
\Delta_r(x) = \{y\in\R^n\colon |y-x|<r\}
\end{equation*}
represent the open ball in $\R^n$ of radius $r$ and center $x$.

\begin{fact} For any element $\sigma$ of $\onrn$, any point $x$ of $\R^n$, and
any positive number $r$, we have $\sigma(\Delta_r(x)) = \Delta_r(\sigma(x))$.
\label{ballfact}
\end{fact}

\noindent This follows directly from the fact that the elements $\sigma$ of
$\onrn$ correspond to rigid motions (isometries) of $\R^n$, i.e.,
$|\sigma(y)-\sigma(x)| = |y-x|$ for any points $x,y\in\R^n$.

\begin{fact} Each element $\sigma$ of $\onrn$ is a homeomorphism of $\R^n$ onto
itself; in particular, $\sigma^{-1}$ is well-defined.
\label{inversefact}
\end{fact}

\noindent Certainly $\sigma$, being an isometry, is continuous. Moreover, it is
easy to see that if $\sigma=(\theta,\xi)$, then
$\tau=(\theta^{-1},-\theta^{-1}\xi)$ is an element of $\onrn$ which inverts the
action of $\sigma$ on $\R^n$. Therefore $\sigma$ is continuously invertible as
well, hence a homeomorphism.

\begin{fact} For any element $\sigma$ of $\onrn$ and any subset $A$ of $\R^n$,
we have
$\sigma(\Int(A)) = \Int(\sigma(A))$.
\label{intfact}
\end{fact}

\noindent This is an immediate consequence of the fact that $\sigma$ is a
homeomorphism of $\R^n$.

%\begin{fact} The open sets in $\M(\R^n)$ are precisely the sets that can be
%obtained from the collection
%\begin{equation*}
%\Big\{ \big\{ \{\sigma_i\} \in \M(\R^n)\colon \sigma_k\in U \big\} \quad
%(k\ge1,\, U\subset\onrn\hbox{ open}) \Big\}
%\end{equation*} by taking finite intersections and arbitrary unions.
%\label{subbasisfact}
%\end{fact}

\begin{fact} Let $D$ be a subset of $\R^n$, and let $\{x_n\}$ be a sequence of
points of $\R^n$, all but finitely many of which belong to $D$. If $\{x_n\}$
converges to some point
$x$, then $x\in\bar D$.
\label{closurefact}
\end{fact}

\begin{fact} Every closed subset of a compact space is itself compact.
\label{clcocofact}
\end{fact}

\begin{fact} In a compact topological space, every sequence has a convergent
subsequence.
\label{subsequencefact}
\end{fact}

\noindent These three statements are simple consequences of elementary point-set
topology; see for instance Munkres~\cite{munkres}, Sections 2.10, 3.5, and
3.7, respectively.
%\noindent These four statements are simple consequences of elementary point-set
%topology; see for instance Munkres~\cite{munkres}, Sections 2.8, 2.10, 3.5, and
%3.7, respectively.

\begin{fact} If $C$ is a compact subset of $\R^n$, then the space $\M(C)$ is
also compact.
\label{mccompactfact}
\end{fact}

\noindent The orthogonal group $O(n)$ is compact (it is clearly
bounded, since each column is a unit vector in $\R^n$ and hence each entry is at
most 1 in absolute value; and it is closed since it is the preimage of the
identity matrix under the continuous map $\theta\mapsto\theta^T\theta$).  Since
$\M(C) = (O(n)\times C)^\infty$, Fact~\ref{mccompactfact} therefore follows from
Tychonov's theorem that arbitrary products of compact spaces are
compact (see \cite[Section 5.1]{munkres}). The compactness of these
spaces $\M(C)$ is crucial to our proofs of Theorems 1 and~2.

\begin{fact} If $\A = \{A_1,A_2,\dots\}$ is a collection of subsets of $\R^n$,
each containing the origin, then $\P(\A,C)$ is a subset of $\M(C)$.
\label{pacinmcfact}
\end{fact}

\noindent We can justify this fact as follows: if $0\in A$ and
$\sigma=(\theta,\xi)$, then $\xi=\xi+\theta(0)\in\sigma(A)$. Thus if
$\sigma(A)\subset C$, we must have $\xi\in C$. Fact~\ref{pacinmcfact} then
follows from the definition~(\ref{pacspace}) of $\P(\A,C)$ by applying this
reasoning to each image $\sigma_i(A_i)$.

\begin{fact} If $\A=\{A_1,A_2,\dots\}$ and ${\cal C}=\{C_1,C_2,\dots\}$ are
collections of subsets of $\R^n$, then $\P(\A,\bigcap_{k=1}^\infty C_k) =
\bigcap_{k=1}^\infty
\P(\A,C_k)$.
\label{paccapfact}
\end{fact}

\noindent This follows immediately from unfolding the definitions of
$\P(\A,\bigcap_{k=1}^\infty C_k)$ and $\bigcap_{k=1}^\infty
\P(\A,C_k)$ using equation~(\ref{pacspace}). In words, Fact~\ref{paccapfact}
states that any packing of $\A$ into the set $\bigcap_{k=1}^\infty C_k$ is
simultaneously a packing of $\A$ into each set $C_k$.

\section{Proofs of Theorems 1 and 2}

\noindent In this section we state the following crucial proposition from which
we deduce Theorems~1 and~2:

\newcommand{\closedproptext}{Let $C$ be a closed subset of $\R^n$, and let
$\A$ be any collection of subsets of $\R^n$. Then the space $\P(\A,C)$ is a
closed subset of $\M(\R^n)$.}
\begin{proposition}
\closedproptext
\label{closedprop}
\end{proposition}

\noindent The proof of Proposition~\ref{closedprop}, while not tricky, is
somewhat long-winded, and therefore we defer it to the next section.
Assuming the validity of Proposition~\ref{closedprop}, we can establish
Theorems~1 and~2 by means of the following lemma:

\begin{lemma} Let $\A=\{A_1,A_2,\dots\}$ and ${\cal C}=\{C_1,C_2,\dots\}$
be collections of subsets of $\R^n$. For each $k\ge1$ define $D_k =
\bigcup_{j=k}^\infty C_j$, and suppose that $D_1$ is bounded. If there exist
packings of $\A$ into $C_j$ for each $j\ge1$, then there exists a packing of $\A$
into the set $\bigcap_{k=1}^\infty \bar D_k$.
\label{alllem}
\end{lemma}

\noindent The set $\bigcap_{k=1}^\infty \bar D_k$ can be compared to the related
set $\bigcap_{k=1}^\infty D_k$, which is simply the lim sup of the sets $C_j$
(the set of all points that are contained in infinitely many of the $C_j$). In
fact, $\bigcap_{k=1}^\infty \bar D_k$ is precisely the set of all points
$x\in\R^n$ such that every neighborhood of $x$ intersects infinitely many of the
$C_j$.

\bigskip
\begin{proof}
By translating the sets $A_i$ if necessary, we may assume that each $A_i$
contains the origin. By hypothesis, there exists a packing of $\A$ into each
$C_j$, so we may choose
\begin{equation*}
T\j \in \P(\A,C_j) \subset \P(\A,\bar D_j) \subset \P(\A,\bar D_1)
\end{equation*} for each $j\ge1$. The set $\bar D_1$ is closed and bounded, hence
compact, and so by Fact~\ref{mccompactfact} the space $\M(\bar D_1)$ is also
compact. Since the sets $A_i$ all contain the origin, the space $\P(\A,\bar
D_1)$ is contained in $\M(\bar D_1)$ by Fact~\ref{pacinmcfact}; we know by
Proposition~\ref{closedprop} that $\P(\A,\bar D_1)$ is a closed set, and so it
is itself compact by Fact~\ref{clcocofact}. Therefore by
Fact~\ref{subsequencefact}, the sequence $\{T_j\}$ of points in $\P(\A,\bar
D_1)$ has a convergent subsequence. By replacing the sequence $\{T_j\}$ by this
subsequence, we may assume that the $T_j$ converge to some element
$T\in\P(\A,\bar D_1)$.

It remains to show that this element $T$ in fact represents a packing of $\A$
into
$\bigcap_{k=1}^\infty \bar D_k$. For each $k\ge1$, the sequence $T\j$ is
contained (except for at most the first $k-1$ terms) in $\P(\A,\bar D_k)$. Since
this set is closed by Proposition~\ref{closedprop}, we see by
Fact~\ref{closurefact} that the limit $T$ is itself an element of $\P(\A,\bar
D_k)$. Because this is true for all $k\ge1$, Fact~\ref{paccapfact} implies
\begin{equation*} T \in \bigcap_{k=1}^\infty \P(\A,\bar D_k) = \P\bigg(
\A,\bigcap_{k=1}^\infty \bar D_k \bigg),
\end{equation*} which establishes the lemma.
\end{proof}

\begin{pflike}{Proof of Theorem~\ref{homothetthm}:} Since $C$ is compact, it is
contained in some ball of radius $R$ centered at the origin, and
therefore each set $(1+\frac1j)C$ is contained in the
ball of radius $2R$ around the origin. Therefore under the hypothesis that there
exist packings of $\A$ into each set $(1+\frac1j)C$, we may apply
Lemma~\ref{alllem} to conclude that there exists a packing of $\A$ into the set
$\bigcap_{k=1}^\infty \bar D_k$, where we have put
\begin{equation}
D_k = \bigcup_{j=k}^\infty (1+{\textstyle\frac1j})C.  \label{Dkdef}
\end{equation}
All that remains to establish the theorem is to show that $\bigcap_{k=1}^\infty
\bar D_k$ is contained in $C$; in other words, we need to show that for every
$x\notin C$, there exists some $k\ge1$ such that $x\notin\bar D_k$.

If $x\notin C$ then, since $C$ is compact (hence closed), there exists a positive
number $\ep$ such that $\Delta_\ep(x)\cap C=\emptyset$. We claim that
\begin{equation}
\hbox{for every }j>2|x|\ep^{-1},\quad \Delta_{\ep/2}(x) \cap
(1+{\textstyle\frac1j})C = \emptyset.
\label{inter}
\end{equation} To see this, suppose that there did exist a point $y$ in
$\Delta_{\ep/2}(x) \cap (1+\frac1j)C$. Since $y\in(1+\frac1j)C$, if we set
$z=(1+\frac1j)^{-1}y$ then
$z\in C$, and by our choice of $\ep$ we therefore have $|x-z|\ge\ep$. On the
other hand, since $y\in\Delta_{\ep/2}(x)$,
\begin{equation*} |x-z| \le |x-y| + |y-z| < \frac\ep2 + |y - (1 +
{\textstyle\frac1j} )^{-1} y | = \frac\ep2 + {|y|\over j+1}.
\end{equation*} The fact that $y\in\Delta_{\ep/2}(x)$ forces $|y|<|x|+\ep/2$,
and so
\begin{equation*} |x-z| < \frac\ep2 + {|x|+\ep/2\over j+1} < \frac\ep2 +
{|x|+\ep/2\over 2|x|/\ep+1} = \ep
\end{equation*} by our choice of $j$. This contradiction establishes
equation~(\ref{inter}).

If we set $k=\lfloor2|x|\ep^{-1}\rfloor+1$, we see from equation~(\ref{inter})
and the definition~(\ref{Dkdef}) of $D_k$ that $\Delta_{\ep/2}(x) \cap D_k =
\emptyset$, which implies that $x\notin\bar D_k$ as desired. This establishes the
theorem.\qed
\end{pflike}

\begin{pflike}{Proof of Theorem~\ref{brickthm}:}
First we make some reductions in the problem. By translating each set $A_i$ if
necessary we may assume that each $A_i$ contains the origin. Similarly, by
translating each brick $B_j$ if necessary, we may assume that each $B_j$ is
contained in the positive orthant of $\R^n$ and has one vertex at the origin,
that is, $B_j = [0,b_{j1}] \times \dots \times [0,b_{jn}]$. Next, by passing to
a suitable subsequence of the $B_j$, we may also assume that $\vol B_j$
decreases monotonically to $V$. At this point we make the assumption that the
dimensions $b_{jm}$ of the bricks $B_j$ are bounded uniformly in $j$ and $m$; at
the end of the proof we will show why this assumption is legitimate. By passing
once again to a suitable subsequence of the $B_j$, we may therefore assume that
for each $1\le m\le n$ the sequence $\{b_{jm}\}$ converges to some number $b_m$,
say.

Since the $b_{jm}$ are uniformly bounded, the sets $B_j$ are all contained in a
single bounded region of $\R^n$, and thus we may apply Lemma~\ref{alllem} to
conclude that there exists a packing of the set $\A$ into $\bigcap_{k=1}^\infty
\bar D_k$, where we have put $D_k = \bigcup_{j=k}^\infty B_j$. The theorem will
therefore be established if we can demonstrate that the intersection
$\bigcap_{k=1}^\infty
\bar D_k$ is contained in the brick $B = [0,b_1] \times \dots \times [0,b_n]$.
For any natural numbers $k$ and $m$ with $1\le m\le n$, define $d_{km} =
\sup_{j\ge k}\{b_{jm}\}$. Then for $j\ge k$ it is clear that $B_j$ is contained
in the closed set $[0,d_{k1}] \times \dots \times [0,d_{kn}]$, and so $\bar D_k$
is contained in the same closed set. Consequently,
\begin{equation*}
\begin{split}
\bigcap_{k=1}^\infty \bar D_k &\subset \bigcap_{k=1}^\infty \big( [0,d_{k1}]
\times \dots \times [0,d_{kn}] \big) \\
&= \big[0,\inf\nolimits_k\{d_{k1}\}\big] \times \dots \times
\big[0,\inf\nolimits_k\{d_{kn}\}\big] \\
&= \big[0,\limsup\nolimits_j\{b_{j1}\}\big] \times \dots \times
\big[0,\limsup\nolimits_j\{b_{jn}\}\big] \\
&= \hskip1pt[0,b_1] \times \dots \times [0,b_n] = B.
\end{split}
\end{equation*}

This establishes the theorem, modulo the assumption that the $b_{jm}$ are
uniformly bounded. This assumption does not hold for a general collection of
bricks of bounded volume, as the simple example $[0,n]\times[0,1/n]$
in $\R^2$ demonstrates. However, in the most natural case---where at least one of
the sets $A_i$ has nonempty interior---we will be able to deduce from the
existence of a packing of $\A$ into each brick $B_j$ that the $b_{jm}$ are
uniformly bounded. In the contrary (less interesting) case, it will also be
possible to reduce to the situation where the $b_{jm}$ are uniformly bounded by
a somewhat different method.

{\narrower
\medskip\noindent {\it Case 1.} At least one of the sets $A_i$ has nonempty
interior.\smallskip

Choose an integer $k$ such that the set $A_k$ has nonempty interior, and then
choose $\eta>0$ such that $A_k$ contains some open ball of radius $\eta$. Since
there exists a packing of $\A$ into each brick $B_j$, we see in particular that
each $B_j$ contains some open ball of radius $\eta$. Certainly then the
dimensions $b_{j1},\dots,b_{jn}$ of each brick $B_j$ must satisfy
$b_{jm}\ge\eta$ for each $1\le m\le n$, and so for each $j\ge1$ and $1\le m\le
n$,
\begin{equation*} 0 < b_{jm} = {\vol B_j\over b_{j1}\dots
b_{j,m-1}b_{j,m+1}\dots b_{jn}}
\le {\vol B_1\over\eta^{n-1}},
\end{equation*} since we have reduced to the case where the $\vol B_j$ are
monotonically decreasing. This shows that the $b_{jm}$ %and hence the $\diam B_j$
are indeed uniformly bounded.

\medskip\noindent {\it Case 2.} All of the $A_i$ have empty interiors.\smallskip

We claim that if there exists a packing of $\A$ into each brick $B_j =
[0,b_{j1}] \times \dots \times [0,b_{jn}]$, then there also exists a packing of
$\A$ into the smaller brick $B'_j = [0,b'_{j1}] \times \dots \times
[0,b'_{jn}]$ where we have defined $b'_{jm} = \min\{b_{jm},\diam B_1\}$. If we
can justify this assertion, the theorem is established in this case as well
since the $b'_{jm}$ are certainly uniformly bounded by $\diam B_1$.

For a collection $\A$ of sets with empty interiors, the packing condition that
the positionings of the sets $A_i$ must have disjoint interiors is no condition
at all; in other words, there exists a packing of the entire collection $\A$ into
$C$ if and only if there exists individual positionings of each set $A_i$ into
$C$. Moreover, we can modify any positioning $\sigma_i(A_i)$ into the brick
$B_j$ so that it becomes a positioning of $A_i$ into $B'_j$, by
taking the rotated/reflected set $\theta_i(A_i)$ and translating
it just enough to lie the positive orthant of~$\R^n$.
More precisely, if $\sigma_i=(\theta_i,\xi_i)$ is such that $\sigma_i(A_i)\subset
B_j$, then we define $\sigma'_i=(\theta_i,\xi'_i)$ where the $m$th coordinate
$\xi'_{im}$ of the vector $\xi'_i\in\R^n$ is given by
\begin{equation*}
\xi'_{im} = \big|\inf \{t\in\pi_i(\theta_i(A_i))\} \big|;
\end{equation*}
here $\pi_i$ denotes the projection map in the $i$th coordinate from $\R^n$ to
$\R$. 

The fact that $\sigma'_i(A_i)$ is contained in the positive orthant of
$\R^n$ follows immediately from the definition of the $\xi'_{im}$. Also, we are
assuming that $A_i$ contains the origin, and so $\xi_i$ is an element of
$\sigma_i(A_i)$; since $\sigma_i(A_i)$ is contained in the positive orthant, it
follows that $\xi'_{im}\le\xi_{im}$, and consequently $\sigma'_i(A_i)$ is
contained in the brick $B_j$. Finally, since $A_i$ contains the origin it is
clear that $\xi'_{im}\le\diam A_i$, and since there exists a packing of $\A$
into $B_1$ we certainly have $\diam A_i\le\diam B_1$. Therefore
$\sigma'_i(A_i)$ is indeed contained in the brick $B'_j$.

Making this modification for each set $A_i$ results in a packing of the entire
collection $\A$ into the smaller brick $B'_j$ (again, the assumption that the
$A_i$ have empty interiors means that we do not need to worry about the relative
positionings of the various $A_i$). As remarked earlier, this justifies the
assumption that the dimensions of our bricks are uniformly bounded, since we may
replace $B_j$ by $B'_j$ throughout.

} % end \narrower

\medskip This completes the proof of the theorem.
\qed
\end{pflike}

In summary, we have established Theorems~1 and~2 modulo a proof of
Proposition~\ref{closedprop}; this proof will be the subject of the following
section.

\section{Proof of Proposition~\ref{closedprop}}

\noindent Proposition~\ref{closedprop} is essentially a consequence of the fact
that the action on $\R^n$ of the space of rigid motions $\onrn$ is continuous.
The following two lemmas, which give concrete statements of the continuity of
this action, will enable us to establish Proposition~\ref{closedprop}. We note
that the space $\onrn$ can in fact be regarded as a metric space, inheriting as
it does the standard metric from $\R^{n^2}\times\R^n$: if $\sigma=(\theta,\xi)$
and $\sigma'=(\theta',\xi')$ are two elements of $\onrn$, then the distance
between them~is
\begin{equation}
d(\sigma',\sigma) = \big( |\theta'-\theta|^2 + |\xi'-\xi|^2 \big)^{1/2} = \bigg(
\sum_{l=1}^n \sum_{m=1}^n (\theta_{lm}'-\theta_{lm})^2 + \sum_{m=1}^n
(\xi_m'-\xi_m)^2 \bigg)^{1/2},
\label{ddef}
\end{equation}
considering $\theta$ and $\theta'$ here simply as $n^2$-tuples of real
numbers rather than elements of $O(n)$.

\begin{lemma}
Let $y$ be a point in $\R^n$ and $U$ be an open subset of $\R^n$. Suppose that
$\sigma$ is an element of $\onrn$ such that $\sigma(y)\in U$. Then there
exists a positive real number $\delta$ such that, for every $\sigma'\in\onrn$
satisfying $d(\sigma',\sigma)<\delta$, we have $\sigma'(y)\in U$.
\label{pointopenlem}
\end{lemma}

\begin{proof} For any $y\in\R^n$ and any pair $\tau=(\theta,\xi)$,
$\tau'=(\theta',\xi')$ of elements of $\onrn$, we have
\begin{equation} | \tau(y) - \tau'(y) | = | \xi + \theta y - \xi' - \theta'y |
\le |
\xi-\xi' | + | (\theta-\theta')y |,
\label{foost}
\end{equation}
We certainly have $|\xi-\xi'| \le
d(\tau,\tau')$ by the definition~(\ref{ddef}) of the metric $d$. On the other
hand, all entries of the matrix $\theta-\theta'$ are also at most
$d(\tau,\tau')$ in absolute value, while the entries of the vector $y$ are at
most $|y|$ in absolute value. Therefore each entry of $(\theta-\theta')y$ is
bounded by $n|y| d(\tau,\tau')$ in absolute value, and so the
inequality~(\ref{foost}) becomes the upper bound
\begin{equation}
| \tau(y) - \tau'(y) | \le d(\tau,\tau') + \bigg( \sum_{m=1}^n \big ( n|y|
d(\tau,\tau') \big)^2 \bigg)^{1/2} = (n^{3/2}|y|+1) d(\tau,\tau')
\label{continuity}
\end{equation}
(we have made no effort to obtain a strong constant in the inequality).

Now if $\sigma$ is an element of $\onrn$ such that $\sigma(y)$ lies in the open
set $U$, then there exists some positive number $\ep$ such that
$\Delta_\ep(\sigma(y)) \subset U$. If we set
$\delta=\ep(n^{3/2}|y|+1)^{-1}$, then for any
$\sigma'\in\onrn$ such that $d(\sigma',\sigma)<\delta$, the upper
bound~(\ref{continuity}) tells us that
\begin{equation*}
|\sigma'(y)-\sigma(y)| \le (n^{3/2}|y|+1) d(\sigma',\sigma) < \ep,
\end{equation*}
and therefore $\sigma'(y) \in \Delta_\ep(\sigma(y)) \subset U$ as desired.
\end{proof}

\begin{lemma}
Let $U_1$ and $U_2$ be open subsets of $\R^n$. Suppose that $\sigma_1$ and
$\sigma_2$ are elements of $\onrn$ such that $\sigma_1(U_1) \cap \sigma_2(U_2)
\ne \emptyset$. Then there exists a positive real number $\delta$ such that, for
every $\sigma_1',\sigma_2'\in\onrn$ satisfying
$d(\sigma_1',\sigma_1)<\delta$ and $d(\sigma_2',\sigma_2)<\delta$, we have
$\sigma_1'(U_1) \cap \sigma_2'(U_2) \ne \emptyset$.
\label{openopenlem}
\end{lemma}

\begin{proof}
Since $\sigma_1(U_1)$ and $\sigma_2(U_2)$ are open sets that are not disjoint, we
can choose a point $x\in\R^n$ and a positive number $\ep$ such that
$\Delta_\ep(x) \subset \sigma_1(U_1) \cap \sigma_2(U_2)$. Using
Fact~\ref{inversefact} we may set $y_1=\sigma_1^{-1}(x)$ and 
$y_2=\sigma_2^{-1}(x)$, so that $\Delta_\ep(y_1)\subset U_1$ and
$\Delta_\ep(y_2)\subset U_2$; we also set
\begin{equation*}
\delta = {\ep \over n^{3/2}\max\{|y_1|,|y_2|\}+1},
\end{equation*}
Then for $i=1$ or 2, for any $\sigma_i'\in\onrn$ such that
$d(\sigma_i',\sigma_i)<\delta$ the upper bound~(\ref{continuity}) tells us that
\begin{equation*}
|\sigma_i'(y_i)-x| = |\sigma_i'(y_i)-\sigma_i(y_i)| \le (n^{3/2}|y_i|+1)
d(\sigma_i',\sigma_i) < \ep,
\end{equation*}
so that $x \in \Delta_\ep(\sigma_i'(y_i)) = \sigma_i'(\Delta_\ep(y_i)) \subset
\sigma_i'(U_i)$ by Fact~\ref{ballfact}. In particular, this shows that $x$ is an
element of $\sigma_1'(U_1) \cap \sigma_2'(U_2)$, which is therefore nonempty as
desired.
\end{proof}

%\medskip\noindent{\bf Proposition \ref{closedprop}.} {\it
%\closedproptext}\medskip

\begin{pflike}{Proof of Proposition~\ref{closedprop}:}
Let $T=\{\tau_i\}$ be a point in $\M(\R^n)\setminus\P(\A,C)$. From the
definition~(\ref{pacspace}) of $\P(\A,C)$, one of the following two cases must
hold.

{\narrower
\medskip\noindent {\it Case 1.} There exists a $k\ge1$ such that
$\tau_k(A_k)\not\subset C$.\smallskip

Choose a point $x\in\tau_k(A_k)\setminus C$, and set $y=\tau_k^{-1}(x)\in A_k$ 
(using Fact~\ref{inversefact}). Applying Lemma~\ref{pointopenlem} with
$\sigma=\tau_k$ and
$U=\R^n\setminus C$, we see that there exists a positive number
$\delta$ such that, for every $\sigma'\in\onrn$ satisfying
$d(\sigma',\tau_k)<\delta$, we have $\sigma'(y)\in \R^n\setminus C$, that is,
$\sigma'(y)\notin C$.

Now define the open neighborhood $\S$ of $T$ in $\M(\R^n)$ by
\begin{equation*}
\S = \big\{ S=\{\sigma_i\}\in\M(\R^n) \colon d(\sigma_k,\tau_k)<\delta \big\}.
\end{equation*}
For every $S\in\S$, we see that $\sigma_k(y)\notin C$ by our choice of $\delta$.
On the other hand, certainly $\sigma_k(y)\in\sigma_k(A_k)$, and so $S$ is not a
packing of $\A$ into $C$. Since this is true for any $S\in\S$, we see that $\S
\subset \M(\R^n) \setminus \P(\A,C)$.

\medskip\noindent {\it Case 2.} There exist positive integers $k\ne l$ such that
$\Int(\tau_k(A_k)) \cap \Int(\tau_l(A_l)) \ne
\emptyset$.\smallskip

Applying Lemma~\ref{openopenlem} with $\sigma_1=\tau_k$, $\sigma_2=\tau_l$,
$U_1=\Int(A_k)$, and $U_2=\Int(A_l)$, we see that there exists a positive real
number $\delta$ such that, for every $\sigma_1',\sigma_2'\in\onrn$ satisfying
$d(\sigma_1',\tau_k)<\delta$ and $d(\sigma_2',\tau_l)<\delta$, we have
\begin{equation*}
\Int(\sigma_1'(A_k)) \cap \Int(\sigma_2'(A_l)) = \sigma_1'(\Int(A_k)) \cap
\sigma_2'(\Int(A_l)) \ne \emptyset
\end{equation*}
(here we have used Fact~\ref{intfact}). Now define the open neighborhood $\S$ of
$T$ in $\M(\R^n)$ by
\begin{equation*}
\S = \big\{ S=\{\sigma_i\}\in\M(\R^n) \colon d(\sigma_k,\tau_k)<\delta \hbox{
and } d(\sigma_l,\tau_l)<\delta \big\}.
\end{equation*}
For every $S\in\S$, we see that $\Int(\sigma_k(A_k)) \cap \Int(\sigma_l(A_l)) \ne
\emptyset$ by our choice of $\delta$, and so $S$ is not a packing of $\A$ with
disjoint interiors. Since this is true for any $S\in\S$, we see that $\S \subset
\M(\R^n) \setminus \P(\A,C)$.

}
\medskip In either case we see that $\M(\R^n)\setminus\P(\A,C)$ contains an open
neighborhood $\S$ of $T$, which shows that $\M(\R^n)\setminus\P(\A,C)$ is an open
set, i.e., $\P(\A,C)$ is a closed subset of $\M(\R^n)$.
\qed
\end{pflike}

\section{Generalizations of Theorems 1 and 2}

\noindent We end by briefly discussing some extensions of Theorems~1
and~2 that can be established by the methods of this paper. First, in
the statements of these two theorems we have claimed that ``packings''
may be replaced by ``oriented packings''. This is true because the
positionings allowed in oriented packings (translations and rotations,
but not reflections) are parametrized by $O(n)^+\times\R^n$, where
$O(n)^+$ is the index-2 subgroup of $O(n)$ consisting of the
orthogonal matrices of determinant 1. Because this subgroup $O(n)^+$
is a compact space in its own right, the analogous statement to
Fact~\ref{mccompactfact} for $\M^+(C) = (O(n)^+\times C)^\infty$ is also true,
and thus all of the arguments of this paper go through for oriented packings upon
simply replacing $\M(C)$ by $\M^+(C)$ at each occurrence.  In the case of
translated packings, where neither rotations nor reflections are allowed, we can
similarly replace each occurrence of $\M(C)$ by $C^\infty$ and the arguments
proceed unchanged (if we like, we can think of the space $C^\infty$ as
$(\{I_n\} \times C)^\infty$, where $\{I_n\}$ is the compact subgroup
of $O(n)$ consisting only of the identity matrix).

It is clear that many variations on Theorems~1 and~2 could be stated
by changing the sequence of sets into which $\A$ can be
packed. The important thing is for this sequence $C_j$ (which is a
shrinking sequence of homothets in Theorem~\ref{homothetthm}, and a
sequence of bricks of varying dimensions in Theorem~\ref{brickthm}) to
have enough structure for the limiting set $\bigcap_{k=1}^\infty \bar
D_k$ to be identified, where $D_k = \bigcup_{j=k}^\infty C_j$ as
defined in the statement of Lemma~\ref{alllem}. This limiting set would
be easy to determine if the $C_j$ were ellipsoids or simplices of
varying dimensions, just to name two possible applications.

Finally we note two ways in which the hypotheses of Theorems~1 and~2
can be weakened. Instead of requiring that the collection $\A$ can be
packed into each set $C_j$, we can require only that for each $j\ge1$
the contracted collection $(1-\frac1j)\A = \{ (1-\frac1j)A_1,
(1-\frac1j)A_2, \dots \}$ can be packed into $C_j$. This is actually
easily seen to be equivalent to the current statements of Theorems~1
and~2. However, we obtain genuinely stronger theorems by weakening the
hypothesis in the following way: for every $j\ge1$, we require only
that the finite collection $\{A_1, \dots, A_j\}$ can be packed into
the set $C_j$. We leave the details of this variation to the reader.

\medskip {\smaller\smaller\baselineskip=12pt
\begin{pflike}{Acknowledgements.}
The author acknowledges the support of Natural Sciences and Engineering Research
Council grant number A5123. The author would also like to thank Mark Hamilton
for his comments on a preliminary version of this paper.
\end{pflike}

%\bigskip
%email .ps to jcta@@math.ucla.edu
%
}

\bibliography{squares}
\bibliographystyle{amsplain}
\end{document}